# New Insights into Integrals Involving Gaussian Sums

Jesus Retamozo M. , Lima-Perú.

**July 22, 2025**

## Abstract

In this article, we explore a series of elementary yet insightful results involving integrals related to Gaussian sums. Using techniques rooted in classical calculus, we derive several identities and evaluate nontrivial definite integrals that emerge naturally in this context. The approach is mostly elementary, avoiding the need for advanced machinery, and aims to shed light on the rich interplay between analysis and arithmetic structures arising from these integrals.

## 1. INTRODUCTION

Within Ramanujan's Collected Papers, two intriguing articles dealing with the calculation of certain integrals, of the type illustrated below:

$$\varphi_w(t) = \int_0^\infty \frac{\cos(\pi t x)}{\cosh(\pi x)} e^{-\pi w x^2} dx$$

$$\theta_w(t) = \int_0^\infty \frac{\cos(\pi t x)}{\cosh(\pi x)} e^{-\pi w x^2} dx$$

In two of his remarkable papers, Ramanujan succeeded in determining the exact value of these integrals for certain values of w and t. However, while he resorted to advanced techniques such as transforms and complex analysis, the present paper seeks to deduce the value of these integrals from elementary calculus, thus offering new intuitions about this class of integrals.

Some of the exact values that Ramanujan obtained are as follows:

$$\int_0^\infty \frac{\cos(2\pi tx)}{\cosh(\pi x)}\cos(\pi x^2)dx = \frac{1+\sqrt{2}\sin(\pi t^2)}{2\sqrt{2}\cosh(\pi t)}$$

$$\int_0^\infty \frac{\cos(2\pi tx)}{\cosh(\pi x)}\sin(\pi x^2)dx = \frac{-1+\sqrt{2}\cos(\pi t^2)}{2\sqrt{2}\cosh(\pi t)}$$

$$\int_0^\infty \frac{\sin(2\pi tx)}{\cosh(\pi x)}\cos(\pi x^2)dx = \frac{\cosh(\pi t)-\cos(\pi t^2)}{2\sinh(\pi t)}$$

$$\int_0^\infty \frac{\sin(2\pi tx)}{\cosh(\pi x)}\sin(\pi x^2)dx = \frac{\sin(\pi t^2)}{2\sinh(\pi t)}$$

## 2. MAIN SECTION

Let's start this section by recalling the following identity, which can be found in the well-known book on identities called Gradshtein and Ryzhik.

$$\int_0^\infty \frac{\cos(px)}{\cosh(x)+\cosh(a)}dx = \frac{\pi\sin(pa)}{\sinh(p\pi)\sinh(a)}$$

We will take advantage of the left-hand side of this equality to deduce some incredible non-trivial integrals and arrive at surprising results.

Making the change of variable $p \to pa$, and then dividing the resulting expression by cosh(a).

$$\int_0^\infty \frac{\cos(pax)}{\cosh(x)+\cosh(a)}dx = \frac{\pi\sin(paa)}{\sinh(pa\pi)\sinh(a)} = \frac{\pi\sin(pa^2)}{\sinh(pa\pi)\sinh(a)}$$

$$\int_0^\infty \frac{\cos(pax)}{\cosh(a)(\cosh(x)+\cosh(a))}dx = \frac{\pi\sin(pa^2)}{\sinh(pa\pi)\cosh(a)\sinh(a)}$$

$$\int_0^\infty \frac{\cos(pax)}{\cosh(a)(\cosh(x)+\cosh(a))}dx = \frac{2\pi\sin(pa^2)}{\sinh(pa\pi)\sinh(2a)}$$

Integrating between infinity and zero with respect to 'a'.

$$\int_0^\infty \int_0^\infty \frac{\cos(pax)}{\cosh(a)(\cosh(x)+\cosh(a))}dx\, da = \int_0^\infty \frac{2\pi\sin(pa^2)}{\sinh(pa\pi)\sinh(2a)}da$$

But on the left side we change the order of integration, and then taking advantage of the symmetry of the integrand, we have:

$$2\int_0^\infty \int_0^\infty \frac{\cos(pax)}{\cosh(a)(\cosh(x)+\cosh(a))}dx\, da$$
$$= \int_0^\infty \int_0^\infty \frac{\cos(pax)}{\cosh(a)(\cosh(x)+\cosh(a))}dx\, da$$
$$+ \int_0^\infty \int_0^\infty \frac{\cos(pax)}{\cosh(x)(\cosh(x)+\cosh(a))}dx\, da$$
$$= \int_0^\infty \int_0^\infty \frac{\cos(pax)}{(\cosh(x)+\cosh(a))}\left(\frac{1}{\cosh(a)}+\frac{1}{\cosh(x)}\right)dx\, da$$

$$= \int_0^\infty \int_0^\infty \frac{\cos(pax)}{\cosh(x)\cosh(a)} dx\, da = \frac{\pi}{2} \int_0^\infty \frac{dx}{\cosh(x) \cosh\left(\frac{\pi}{2}px\right)}$$

Replacing in the above.

$$\int_0^\infty \frac{2\pi \sin(pa^2)}{\sinh(pa\pi)\sinh(2a)} da = \frac{\pi}{4} \int_0^\infty \frac{dx}{\cosh(x) \cosh\left(\frac{\pi}{2}px\right)}$$

$$\int_0^\infty \frac{\sin(pa^2)}{\sinh(pa\pi)\sinh(2a)} da = \frac{1}{8} \int_0^\infty \frac{dx}{\cosh(x) \cosh\left(\frac{\pi}{2}px\right)}$$

This is the main identity we needed to arrive at the surprising results, you can see the difference between the methods used by Ramanujan, these are simpler.

## 3. EXACT VALUES

Back to the main identity and making $p = \frac{2}{\pi}$.

$$\int_0^\infty \frac{\sin(\frac{2}{\pi}a^2)}{\sinh(2a)\sinh(2a)} da = \frac{1}{8} \int_0^\infty \frac{dx}{\cosh(x)\cosh(x)}$$

$$\int_0^\infty \frac{\sin(\frac{2}{\pi}a^2)}{\sinh(2a)^2} da = \frac{1}{8} \int_0^\infty \frac{dx}{\cosh(x)^2}$$

The integral of the right side is elementary and so we arrive at:

$$\int_0^\infty \frac{dx}{\cosh(x)^2} = \int_0^\infty \mathrm{sech}(x)^2\, dx = \tanh(x) \bigg|_0^\infty = 1$$

Finally:

$$\int_0^\infty \frac{\sin(\frac{2}{\pi}a^2)}{\sinh(2a)^2} da = \frac{1}{8}$$

With $p = \frac{4}{\pi}$:

$$\int_0^\infty \frac{\sin(\frac{4}{\pi}a^2)}{\sinh(4a)\sinh(2a)} da = \frac{1}{8} \int_0^\infty \frac{dx}{\cosh(x)\cosh(2x)}$$

The right side is elementary but a bit laborious, so omitting that process, we arrive at the following result.

$$\int_0^\infty \frac{\sin(\frac{4}{\pi}a^2)}{\sinh(4a)\sinh(2a)} da = \frac{\pi}{16}(\sqrt{2} - 1)$$

Something interesting happens when we make p equal to zero in the principal formula.

$$\lim_{p \to 0} \int_0^\infty \frac{\sin(pa^2)}{\sinh(pa\pi)\sinh(2a)} da = \frac{1}{8} \int_0^\infty \frac{dx}{\cosh(x) \cosh\left(\frac{\pi}{2}(0)x\right)}$$

$$\int_0^\infty \lim_{p \to 0} \frac{\sin(pa^2)}{\sinh(pa\pi)\sinh(2a)} da = \frac{1}{8} \int_0^\infty \frac{dx}{\cosh(x)}$$

The right-hand side is trivial and on the left-hand side we have that:

$$\frac{1}{\pi} \int_0^\infty \frac{a}{\sinh(2a)} 7 da = \frac{\pi}{16}$$

$$\frac{1}{4} \int_0^\infty \frac{x}{\sinh(x)} dx = \frac{\pi^2}{16} \quad \to \quad \int_0^\infty \frac{x}{\sinh(x)} dx = \frac{\pi^2}{4}$$

Decomposing the left side in infinite series.

$$2 \sum_{k=0}^\infty \int_0^\infty x e^{-x(2k+1)} dx = \frac{\pi^2}{4}$$

$$\sum_{k=0}^\infty \frac{1}{(2k+1)^2} = \frac{\pi^2}{8}$$

We have solved the Basel problem, and we can also deduce a connection between this large problem and the main integral of this article.

## 4. ANOTHER INTRIGUING IDENTITY

Returning to the following identity:

$$\int_0^\infty \frac{\cos(px)}{\cosh(x) + \cosh(a)} dx = \frac{\pi \sin(pa)}{\sinh(p\pi)\sinh(a)}$$

Making the change of variable $p \to pa$ and then integrating between infinity and zero.

$$\int_0^\infty \int_0^\infty \frac{\cos(pax)}{\cosh(x) + \cosh(a)} dx \, da = \int_0^\infty \frac{\pi \sin(pa^2)}{\sinh(pa\pi)\sinh(a)} da$$

On the right hand side, in the denominator we use the summation identity for hyperbolic cosines, and in the numerator we use the following algebraic identity:

$$xa = \frac{1}{4}((x+a)^2 - (x-a)^2)$$

$$\to \int_0^\infty \int_0^\infty \frac{\cos(pax)}{\cosh(x) + \cosh(a)} dx \, da = \frac{1}{4} \int_{-\infty}^\infty \int_{-\infty}^\infty \frac{\cos(pax)}{\cosh(x) + \cosh(a)} dx \, da$$

$$= \frac{1}{4} \int_{-\infty}^\infty \int_{-\infty}^\infty \frac{\cos\left(\frac{p}{4}((x+a)^2 - (x-a)^2)\right)}{2 \cosh\left(\frac{x+a}{2}\right) \cosh\left(\frac{x-a}{2}\right)} dx \, da$$

$$= \frac{1}{2} \int_{-\infty}^\infty \int_{-\infty}^\infty \frac{\cos(p((x+a)^2 - (x-a)^2))}{\cosh(x+a) \cosh(x-a)} dx \, da$$

Using variable changes:

$$y = x + a, z = x - a$$

$$= \frac{1}{2}\int_{-\infty}^{\infty}\int_{-\infty}^{\infty}\frac{\cos(p(y^2-z^2))}{\cosh(y)\cosh(z)}\left(\frac{dydz}{2}\right) = \int_{0}^{\infty}\int_{0}^{\infty}\frac{\cos(p(y^2-z^2))}{\cosh(y)\cosh(z)}dydz$$

$$= \int_{0}^{\infty}\int_{0}^{\infty}\frac{\cos(py^2)\cos(pz^2) + \sin(py^2)\sin(pz^2)}{\cosh(y)\cosh(z)}dydz$$

$$= \int_{0}^{\infty}\int_{0}^{\infty}\frac{\cos(py^2)\cos(pz^2)}{\cosh(y)\cosh(z)}dydz + \int_{0}^{\infty}\int_{0}^{\infty}\frac{\sin(py^2)\sin(pz^2)}{\cosh(y)\cosh(z)}dydz$$

$$= \left(\int_{0}^{\infty}\frac{\cos(py^2)}{\cosh(y)}dy\right)^2 + \left(\int_{0}^{\infty}\frac{\sin(py^2)}{\cosh(y)}dy\right)^2$$

Finally:

$$\int_{0}^{\infty}\frac{\pi\sin(pa^2)}{\sinh(pa\pi)\sinh(a)}da = \left(\int_{0}^{\infty}\frac{\cos(py^2)}{\cosh(y)}dy\right)^2 + \left(\int_{0}^{\infty}\frac{\sin(py^2)}{\cosh(y)}dy\right)^2$$

Again, making p=0 we can solve again the Basel problem.

$$\int_{0}^{\infty}\frac{a}{\sinh(a)}da = \left(\int_{0}^{\infty}\frac{1}{\cosh(y)}dy\right)^2 + (0)^2 = \frac{\pi^2}{4}$$

In the new main identity, for certain values of 'p' we can deduce the exact value of the integral on the left side, we can see how a simple change in comparison with the previous main result changed very drastically the right side, while the first one seemed to yield exact values easily, the second one requires more work but for certain values of 'p' it can be deduced exactly.

## 5. CONCLUSIONS

In the present work, we have demonstrated the possibility of deducing classical results, previously established by Srinivasa Ramanujan, using only elementary calculus techniques. This approach has not only allowed us to revisit these fundamental truths from a simplified perspective, but has also led to the discovery and deduction of new identities, thus enriching the existing body of knowledge.

The elegance and inherent 'magic' revealed in the two identities deduced in this article suggest a vast potential for future exploration, with the promise of leading to even more interesting results. However, to address the generalization of these findings and explore their full scope, it is postulated that the use of more advanced tools, such as complex variable analysis, will be indispensable.